\documentclass[a4paper]{article}
\usepackage[utf8]{inputenc}
\usepackage{mathtools}
\usepackage{amsthm}
\usepackage{amssymb}
\usepackage{parskip}
\usepackage{natbib}
\usepackage{hyperref}
\usepackage{xcolor}
\usepackage[a4paper, margin=2.5cm]{geometry}
\usepackage{tikz-cd}
\usepackage{authblk}

\newtheorem{theorem}{Theorem}
\newtheorem{lemma}[theorem]{Lemma}
\newtheorem{remark}[theorem]{Remark}
\newtheorem{definition}[theorem]{Definition}
\newtheorem{corollary}[theorem]{Corollary}

\title{Duality for Neural Networks through Reproducing Kernel Banach Spaces}
\author[1,*]{Len Spek}
\author[1]{Tjeerd Jan Heeringa}
\author[2]{Felix Schwenninger}
\author[1]{Christoph Brune}
\affil[1]{Mathematics of Imaging \& AI, University of Twente, Enschede, The Netherlands}
\affil[2]{Mathematics of Systems Theory, University of Twente, Enschede, The Netherlands}
\date{February 2023}

\begin{document}
\maketitle

\begin{abstract}
    \noindent Reproducing Kernel Hilbert spaces (RKHS) have been a very successful tool in various areas of machine learning. Recently, Barron spaces have been used to prove bounds on the generalisation error for neural networks. Unfortunately, Barron spaces cannot be understood in terms of RKHS due to the strong nonlinear coupling of the weights. This can be solved by using the more general Reproducing Kernel Banach spaces (RKBS). We show that these Barron spaces belong to a class of integral RKBS. This class can also be understood as an infinite union of RKHS spaces. Furthermore, we show that the dual space of such RKBSs, is again an RKBS where the roles of the data and parameters are interchanged, forming an adjoint pair of RKBSs including a reproducing kernel. This allows us to construct the saddle point problem for neural networks, which can be used in the whole field of primal-dual optimisation.
    \end{abstract}

\section{Introduction}
Neural networks are often considered to be black boxes from a mathematical perspective. However, a lot of recent progress has been made by \citet{e_priori_2019} by proving generalisation bounds and approximation errors for two layer neural networks. The key tool used in these proofs is the Barron space, a space of functions which can be represented by an infinite width neural networks with bounded weights. 

Reproducing Kernel Hilbert spaces (RKHS) have been a very useful tool in machine learning methods. The kernel property, together with the easy representation of the dual space due to the inner product structure, give easy access to various results from representation and optimisation theory. For example, this is the key to prove the representer theorem for random feature models \citep{bach_equivalence_2017}. However, to consider full two layer neural networks \citet{e_barron_2021} has shown that one needs an infinite union of RKHS.

In this paper, we show that the Barron spaces satisfy the Reproducing Kernel Banach space (RKBS) property. These spaces were introduced by \citet{zhang_reproducing_2009} to study the representation properties in machine learning, and were brought into the context of neural networks by \citet{bartolucci_understanding_2023} to prove a representer theorem. In merging these two approaches, Barron spaces and RKBS, the approximation error bounds of the Barron space and the representation results of RKBS can be combined. As an RKBS is not by definition a Hilbert space, a description of the dual space is necessary to study the eponymous reproducing kernel. 

The class of RKBS used by \citet{bartolucci_understanding_2023} for neural networks have an integral representation with measure over the parameters. We show that its dual space can be identified with an RKBS which has an integral representation with measures over the data. These two spaces form an adjoint RKBS pair with a well-defined reproducing kernel. This enables the fundamental basis for learned primal-dual methods for empirical risk minimisation in the context of architecture search and targeted optimisation. Here the primal variables are the neural networks parameters and the dual variables are observations of the data. The dual view opens a door to the field of primal-dual algorithms. 

As mentioned above, the Barron space can be written as an infinite union of RKHS. We show that this argument generalises to RKBS, where an RKBS with an integral representation over measures can be written as an infinite union of certain spaces with an $L^p$-norm, which also satisfy the RKBS property. Surprisingly, this result is independent of the choice of $p$. For $p=2$, these spaces satisfy the RKHS property, and we recover the result of \citet{e_barron_2021}. 

\subsection{Research Context}
The Reproducing Kernel Hilbert spaces have a long history, \citep{aronszajn_theory_1950}, and have been developed into an indispensable concept in several fields of mathematics, such as operator theoretic approaches in complex analysis. In the last decades, they have also become very relevant in machine learning \citep{bach_equivalence_2017}. Advantages of these spaces are that inputs can be explained geometrically and that many powerful results from functional analysis can be used \citep{scholkopf_learning_2018}. This has led to many effective algorithms for different learning tasks, such as support vector machines \citep{wahba_support_1999,evgeniou_regularization_2000}. 

However, the Hilbert structure does not appear to be suitable for certain learning tasks. These spaces offer effectively only a single geometry, as all Hilbert spaces of the same dimension are isometrically isomorphic, cf \cite{rudin_functional_1991}. On the contrary, general Banach spaces offer much more flexibility, as the norm does not need to follow the parallelogram law. This is relevant when considering $l^1$-regularisation, which is useful for creating sparsity \citep{tibshirani_regression_1996,micchelli_feature_2007,song_reproducing_2013}. This naturally leads to the question if an equivalent theory for Reproducing Kernel Banach spaces can be developed. Such a generalisation was proposed by \citet{canu_functional_2003}, where the Reproducing Kernel Banach space is based on bounded point evaluation of functions. I suggest to rather say "to consider semi inner products in Banach spaces, which to some extent allow mimicking the Hilbert space case \citep{lumer_semi-inner-product_1960,giles_classes_1967} and are applied to machine learning \citep{zhang_learning_2017,zhang_reproducing_2009,garcia_sampling_2013}. A good overview on the general RKBS theory can be found in \citet{lin_reproducing_2022}.

Neural networks are often studied in the infinite width case, as the universal approximation theorem guarantees that any continuous function can be approximated in the supremum norm, given some mild assumptions on the activation function \citep{cybenko_approximation_1989,hornik_approximation_1991,leshno_multilayer_1993}. \citet{barron_universal_1993} managed to prove an $L^2$ convergence rate based on Fourier techniques, which do not suffer from the curse of dimensionality. In honour of his work, \citet{e_towards_2020} defined the Barron spaces, which consists of functions which can be approximated by neural networks without the weights blowing up. With the corresponding norms, \citet{e_priori_2019} were able to prove bounds of different kinds of generalisation errors. 

On the other hand, \citet{bartolucci_understanding_2023} used RKBS to understand neural networks and prove a general representer theorem for these spaces. The norm of the RKBS corresponding to ReLU activation functions can be represented in terms of an inverse Radon transform of a bounded real measure \citep{parhi_banach_2021}.

Duality is a key concept in the study of functional analysis. The duality pairing for Banach spaces is the general analogue of the inner product of Hilbert spaces. This enables primal-dual algorithms, which have been very successful in the field of convex analysis. For example, total variation denoising methods have been significantly improved using a proper description of dual space of functions of bounded variation \citep{chambolle_algorithm_2004,chan_nonlinear_1999,osher_iterative_2005, brune_primal_2011}. For the gradient descent algorithm, advances are made to incorporate duality, such dual coordinated ascent methods \citep{shalev-shwartz_stochastic_2013, raj_explicit_2021}.  

\subsection{Overview}
In Section~\ref{sec:RKBS}, we first give an overview of the different, but equivalent definitions of RKHS and how they generalise to RKBS. Next, we introduce the two spaces corresponding to infinitely wide neural networks: The class of integral RKBS from \citet{bartolucci_understanding_2023} and the Barron spaces of \citet{e_priori_2019}. We then show that the Barron space indeed satisfies the RKBS property.

In Section~\ref{sec:union}, we show that an RKBS with an integral representation over measures can be written as an infinite union of certain $L^p$-like spaces, which also satisfy the RKBS property.

In Section~\ref{sec:duality_rkbs}, we prove that an RKBS which has an integral representation with measure over the parameters, has a dual space which is an RKBS with an integral representation with measures over the data. Using this, we construct the saddle point problem for an unconstrained optimisation problem.

\subsection{Notation}
\label{sec:notation}
In the following sections, we employ the following notation conventions and assumptions: $X \subseteq \mathbb{R}^d$ denotes the set from which the data points $x\in X$ are sampled and $\Omega\subseteq \mathbb{R}^D$ denotes the set of weights $w\in \Omega$ over which we want to optimise our model. We assume that $X$ and $\Omega$ are closed and endowed with the natural subspace topology.

We consider a bounded, continuous kernel $\varphi\in C_0(X\times\Omega)$ which vanishes at infinity. We introduce the notation $\varphi_w\in C_0(X)$ for the function $x \mapsto \varphi(x,w)$ given some $w\in \Omega$. Similarly, we write $\varphi_x\in C_0(\Omega)$ for the function $w \mapsto \varphi(x,w)$ given some $x\in X$.

We consider real-valued Radon measures, i.e. the signed, finite, regular Borel measures, defined on either $X$ or $\Omega$. The natural norm for a Radon measure $\mu$ is the total variation norm.
\[\|\mu\|_{\mathcal{M}(\Omega)}:= |\mu|(\Omega)\]
where $|\mu|$ is defined through the Jordan decomposition: There exists two unique positive measures $\mu^+$ and $\mu^-$ such that $\mu = \mu^+ - \mu^-$. Then $|\mu| := \mu^+ + \mu^-$. The Banach space $\mathcal{M}(\Omega)$ is defined as the set of all Radon measures on $\Omega$ with finite total variation norm. A measure $\mu \in \mathcal{M}(\Omega)$ represents the distributions of weights. Similarly, $\rho \in \mathcal{M}(X)$ represents a data distribution. The subspace of probability measures is denoted by $\mathcal{P}(\Omega),\mathcal{P}(X)$ respectively. Point measures will always be referred to by using the Greek letter $\delta$.

The space of continuous functions vanishing at infinity, $C_0(\Omega)$, is a Banach space when endowed with the supremum norm. The Riesz representation theorem implies that $\mathcal{M}(X),\mathcal{M}(\Omega)$ can be identified as the dual spaces of $C_0(X),C_0(\Omega)$ respectively, and the duality pairings imply that the following integrals are well-defined and finite
\begin{align*}
\langle \rho, \varphi_w \rangle &= \int_X \varphi_w(x) d\rho(x)\\
\langle \mu, \varphi_x \rangle &= \int_\Omega \varphi_x(w) d\mu(w)\\
\end{align*}
for all $\rho \in \mathcal{M}(X), \mu \in \mathcal{M}(\Omega)$.

Given some measure $\mu \in \mathcal{M}(\Omega)$ and integrable function $h\in L^1(\mu)$, we use the notation $d\nu := hd\mu$ to define a measure $\nu \in \mathcal{M}(\Omega)$ which is absolutely continuous with respect to $\mu$ and for every Borel set $A\subseteq \Omega$
\[\nu(A) = \int_A h(w) d\mu(w)\]
Finally, when there exists an isometric isomorphism between two Banach spaces $\mathcal{B}_1,\mathcal{B}_2$, we denote this by $\mathcal{B}_1 \cong \mathcal{B}_2$. 

\section{Reproducing Kernel Banach Spaces}
\label{sec:RKBS}
In this section, we will introduce the Reproducing Kernel Banach spaces (RKBS) similarly as \citet{bartolucci_understanding_2023} by first comparing three equivalent definitions of Reproducing Kernel Hilbert spaces (RKHS). We will then introduce a class of integral RKBS, which will be the main object of study in this paper. Finally, we show that the Barron space with ReLU functions of \citet{e_priori_2019} satisfies the integral RKBS properties. 

\subsection{Definition of RKHS}
Reproducing Kernel Hilbert spaces can be equivalently characterised in three ways: boundedness of point evaluation, existence of a reproducing kernel and as a quotient space of a feature space. 

First, the classical definition where point evaluation is a bounded functional.
\begin{definition}\label{def:rkhs1}
A Hilbert space $\mathcal{H}$ of functions on $X$ is a reproducing kernel Hilbert space (RHKS) if, for all $x\in X$, there exists a constant $C_{x}>0$ such that for all $f\in \mathcal{H}$
\begin{equation}
|f(x)|\leq C_x \|f\|_{\mathcal{H}}
\end{equation}
\end{definition}

As point evaluation is a linear functional, the Riesz representation theorem guarantees the existence of the eponymous reproducing kernel 
\begin{theorem}\label{def:rkhs2}
A Hilbert space $\mathcal{H}$ of functions on $X$ satisfies the RKHS property if and only if there exists a function $K: X\times X \mapsto \mathbb{R}$ such that for all $x\in X$
\begin{equation}
\begin{split}
    &K(x,\cdot) \in \mathcal{H}\\
    &f(x) = \langle f, K(x,\cdot)\rangle \qquad \forall f\in \mathcal{H}    
\end{split}
\end{equation}
\end{theorem}
This is easy to check that $K$ is symmetric and positive definite. Each such kernel defines a unique RKHS \citep{aronszajn_theory_1950}. This property of an RKHS has proven very useful in empirical risk minimisation problems, where the minimiser can be written in terms of the reproducing kernel.

Finally, we introduce a definition which a popular in machine learning. An RKHS can be written as a quotient space of a feature space $\Psi$
\begin{theorem}\label{def:rkhs3}
A Hilbert space $\mathcal{H}$ of functions on $X$ satisfies the RKHS property if and only if there exists a Hilbert space $\Psi$ and a map $\psi: X \mapsto \Psi$ such that 
\begin{equation}
\begin{split}
    \mathcal{H} = \Psi/\mathcal{N}(A)\\
    \|f\|_{\mathcal{H}} = \inf_{f=A\nu} \|\nu\|_{\Psi}
\end{split}
\end{equation}
where $A$ maps features in $\Psi$ to functions on $X$ and is defined as
\begin{equation}
    (A\nu)(x) = \langle \psi(x), \nu \rangle
\end{equation}
for all $x\in X$ and $\nu \in \Psi$.
\end{theorem}

The reproducing kernel $K$ is related to the map $\psi$ via the adjoint of the mapping $A:\Psi\mapsto\mathcal{H}$
\begin{equation}
    \psi(x) = A^* K(x,\cdot)
\end{equation}

\subsection{Introduction to RKBS}
Definition~\ref{def:rkhs1} immediately translates to a Banach setting.
\begin{definition}\label{def:rkbs1}
A Banach space $\mathcal{B}$ of functions on $X$ is a reproducing kernel Hilbert space (RHKS) if, for all $x\in X$, there exists a constant $C_{x}>0$ such that for all $f\in \mathcal{B}$
\begin{equation} \label{eq:rkbs_def}
|f(x)|\leq C_x \|f\|_{\mathcal{B}}
\end{equation}
\end{definition}
A trivial example of an RKBS, besides any RKHS space, is the Banach space of bounded continuous functions equipped with to the supremum norm. 

\begin{remark}
Any Banach space $\mathcal{B}$ can be identified with the subspace $j(\mathcal{B})$ of its bi-dual $\mathcal{B}^{**}$ via the canonical embedding $j$. Since $j(\mathcal{B})$ is a space of functions on $\mathcal{B}^*$ and by definition satisfies \eqref{eq:rkbs_def}, it is an RKBS on $\mathcal{B}^*$. For Hilbert spaces, this means that any Hilbert space $\mathcal{H}$ can be identified as an RKHS of functions over $\mathcal{H}$.
\end{remark}

The canonical characterisation is, however, not very useful in practise. As in most applications, $X$ is assumed to be relatively `small' compared to $\mathcal{B}$. One sufficient condition for this intuitive constraint is that $X$ is chosen such that $\mathcal{B}$ contains non-linear functions of $X$. 

Theorem~\ref{def:rkhs3} also translates to the Banach setting, by replacing the inner product with a duality pairing between $\Psi$ and its dual space.
\begin{theorem}\citep[Proposition 3.3]{bartolucci_understanding_2023}\label{def:rkbs3}
A Banach space $\mathcal{B}$ of functions on $X$ satisfies the RKBS property if and only if there exists a Banach space $\Psi$ and a map $\psi: X \mapsto \Psi^*$ such that 
\begin{equation}
\begin{split}
    \mathcal{B} = \Psi/\mathcal{N}(A)\\
    \|f\|_{\mathcal{B}} = \inf_{f=A\nu} \|\nu\|_{\Psi}
\end{split}
\end{equation}
where the linear transformation $A$ maps elements of the Banach space $\Psi$ to functions on $X$ and is defined as
\begin{equation}
    (A\nu)(x) := \langle \psi(x), \nu \rangle
\end{equation}
for all $x\in X$ and $\nu \in \Psi$.
\end{theorem} 
Note that the null-space $\mathcal{N}(A)$ is closed, so the quotient space $\Psi/\mathcal{N}(A)$ is complete and hence a proper Banach space. 

Unfortunately, Definition~\ref{def:rkhs2} does not generalise easily to the Banach setting due to the lack of an analogue to the Riesz-Frechet representation theorem. However, based on Definition~\ref{def:rkbs1}, the point evaluation operator $\delta_x$ is an element of the dual space of $\mathcal{B}$, which allows for a definition of an adjoint RKBS with a reproducing kernel.
\begin{definition}\citep[Definition 2.2]{lin_reproducing_2022}\label{def:rkbs2}
Let $\mathcal{B}$ be an RKBS on the set $X$. If there exists a Banach space $\mathcal{B}^{\#}$ of functions on a set $\Omega$, such that it can be isomorphically embedded into $\mathcal{B}^*$ and if there exists a function $K:X \times \Omega \mapsto \mathbb{R}$, such that $K(x,\cdot)\in \mathcal{B}^{\#}$ for all $x\in X$ and
\begin{equation}
    f(x)= \langle K(x,\cdot), f\rangle 
\end{equation}
for all $x\in X$ and $f\in \mathcal{B}$, then we call $K$ a reproducing kernel for $\mathcal{B}$.

If, in addition, $\mathcal{B}^{\#}$ is also an RKBS on $\Omega$ and it holds that $K(\cdot, w)\in \mathcal{B}$ for all $w\in \Omega$ and
\begin{equation}
    g(w) = \langle g, K(\cdot, w)\rangle 
\end{equation}
for all $w\in \Omega$ and $g\in \mathcal{B}^{\#}$, then we call $\mathcal{B}^{\#}$ an adjoint RKBS of $\mathcal{B}$ and call $\mathcal{B}$ and $\mathcal{B}^\#$ an adjoint pair of RKBS. In this case $K^*(w,x):= K(x,w)$, for $x\in X$ and $w\in \Omega$, is a reproducing kernel of $\mathcal{B}^{\#}$.
\end{definition}
Note that compared to the RKHS case, the reproducing kernel $K$ is no longer symmetric. Also, it follows immediately that any RKHS is adjoint to itself. Furthermore, the adjoint RKBS needs not be unique, as it depends on the chosen set $\Omega$ and the representation of the space $\mathcal{B}^{\#}$. 

For these adjoint pairs of RKBS, we recover again the identity
\begin{equation}
    \psi(x) = A^* K(x,\cdot)
\end{equation}

\citet[Theorem 2.4]{lin_reproducing_2022} prove that any separable RKBS admits a reproducing kernel. However, the separability conditions are not necessary, as there exists a canonical way to define the reproducing kernel for any RKBS by using that the elements of the dual space $\mathcal{B}^*$ can be understood as functions of $\mathcal{B}$.
\begin{theorem}
Let $\mathcal{B}$ be an RKBS on the set $X$ and let $\delta_x\in \mathcal{B}^*$ be the point evaluation functional at $x\in X$. Then 
\begin{equation}
 \mathcal{B}^{\#}:= \overline{\mathrm{span}\{\delta_x|x\in X\}}
\end{equation}
and $\mathcal{B}$ form an adjoint pair of RKBS with a kernel $K: X\times \mathcal{B}\mapsto \mathbb{R}$
\begin{equation}
K(x,f) := f(x)    
\end{equation}
for all $x\in X$, $f\in \mathcal{B}$.
\end{theorem}
\textit{Proof.} By definition, $\mathcal{B}^{\#}$ is a Banach space of functions of $\mathcal{B}$. Furthermore, $K(x,\cdot) = \delta_x \in \mathcal{B}^{\#}$ for all $x\in X$, and hence for all $f\in \mathcal{B}$
\[\langle f, K(x,\cdot) \rangle = \langle f, K(x,\cdot) \rangle = f(x)\]
Moreover, $K(\cdot,f)=f\in \mathcal{B}$ and for hence for all $g\in \mathcal{B}^{\#}$ 
\[\langle g, K(\cdot, f) \rangle = \langle g, f \rangle = g(f)\]
by definition of the duality pairing. \qed

Again, in most applications, we want $\Omega$ to also be 'small' or that $\mathcal{B}^{\#}$ also contains non-linear functions of $\Omega$. In the rest of this paper, $X$ and $\Omega$ are always 'small' with respect to the different RKBSs. 

\subsection{A Class of Integral RKBS} 
Next, we will describe a class of RKBS where we take the Radon measures as the feature space. As we will see in the next section, this is relevant to discuss spaces of functions defined by neural networks. We will construct this RKBS using the feature map of Definition~\ref{def:rkbs3}.

We fix a function $\varphi \in C_0(X\times \Omega)$. For the feature space $\Psi$, we choose the space of Radon measures $\mathcal{M}(\Omega)$ equipped with the total variation norm. Then we define the feature map $\psi(x):= \varphi_x$, where we implicitly assume that $\varphi_x\in C_0(\Omega)$ is canonically embedded into the bi-dual of $ C_0(\Omega)$. 
\begin{equation}\label{eq:A}
    (A\mu)(x) := \langle \psi(x), \mu \rangle =\int_\Omega \varphi_x(w) d\mu(w) = \int_\Omega \varphi(x,w) d\mu(w)
\end{equation}
for all $\mu\in \mathcal{M}(\Omega)$. The RKBS $\mathcal{F}(X,\Omega)$ is then defined as 
\begin{equation}\label{eq:F-space}
\mathcal{F}(X,\Omega):= \{f:X\mapsto \mathbb{R} |\, \exists \mu \in \mathcal{M}(\Omega) \text{ s.t. } f=A\mu \}    
\end{equation}
Hence, $\mathcal{F}(X,\Omega)$ can be identified as the quotient space $\mathcal{M}(\Omega)/\mathcal{N}(A)$ and its norm is given by 
\begin{equation}\label{eq:norm_f}
\|f\|_{\mathcal{F}(X,\Omega)} := \inf_{f=A\mu} \|\mu\|_{\mathcal{M}(\Omega)}
\end{equation}
In plain terms, $\mathcal{F}(X,\Omega)$ contains functions $f$ which can be represented by an integral over $\varphi$ with any measure $\mu$ and the norm of $f$ depends on how large this measure, in terms of its total variation, needs to be. Note that $\varphi \in C_0(X\times \Omega)$ in combination with dominated convergence imply that $\mathcal{F}(X,\Omega)$ is a subspace of $C_0(X)$, albeit with a different norm. 

This integral RKBS $\mathcal{F}(X,\Omega)$ is equivalent to the space $\mathcal{F}_1$ with variation norm $\gamma_1$ used by \citet{bach_breaking_2017} and \citet{chizat_implicit_2020}. The subscript $1$ denotes that the norm behaves similarly to a $1$-norm. Similarly, this space is equivalent to the integral RKBS defined by \citet{bartolucci_understanding_2023}, where $\varphi(x,w)$ is replaced by $\rho(x,\theta)\beta(\theta)$. The smoothing function $\beta$ is used to ensure that the integral \eqref{eq:A} converges for all $\mu$. However, we opted to combine $\rho$ and $\beta$ into $\varphi$, as it produces cleaner notation in later sections. For neural networks with ReLU activation functions, \citet{parhi_banach_2021} defined a norm in terms of the inverse Radon transform of the network, which is equivalent to the $\mathcal{F}(X,\Omega)$-norm \citep{bartolucci_understanding_2023}. We will add to this list of equivalences by proving that the Barron spaces defined by \citet{e_priori_2019} are also an instance of an integral class of RKBS. 

\subsection{Barron Spaces}
The Barron norm for neural networks is defined by \citet{e_priori_2019} to quantify the size of the weights in a neural network of infinite width. In this section we will repeat this definition and show that it satisfies the RKBS property, Definition~\ref{def:rkbs3}, by showing it is equivalent to the RKBS $\mathcal{F}(X,\Omega)$ for a certain $\varphi$.

A shallow neural network or perceptron $f$ maps data vectors $x\in X$, where we assume $X$ is compact, to a scalar output. We consider a continuous activation function $\sigma:\mathbb{R}\mapsto \mathbb{R}$ and some weights $(a_j,v_j,b_j) \in \mathbb{R}\times \Omega\subseteq \mathbb{R}^{d+1}$, where $j\in\{ 1,\cdots m\}$, and functions of the form 
\begin{equation}
    f(x) = \frac{1}{m} \sum_{j=1}^m a_j \sigma(v_j^T x + b_j)
\end{equation}
Formally transitioning to the continuous limit, we can consider functions $f:X\to \mathbb{R}$, for which there exists a probability measure $\pi\in \mathcal{P}(\mathbb{R}\times\Omega)$ such that
\begin{equation}
    f(x) := (A\pi)(x) := \int_{\mathbb{R}\times \Omega} a \sigma(v^T x + b) d\pi(a,v,b)
\end{equation}
The vector space  of such function is called Barron space. The Barron norm of $f$ is defined by taking the infimum of all probability measures $\pi$, which produce $f$, more precisely,
\begin{equation}
\|f\|_{\mathcal{B}_\sigma} := \inf_{f=A\pi} \int_{\mathbb{R}\times \Omega} |a|(1+\|v\|_1+|b|) d\pi(a,v,b)
\end{equation}
For the Barron space to be well-defined, we either require that $\Omega$ is compact or that $\sigma$ grows at most linearly, i.e. $\sigma(y)\in O(y)$. In this case there exists an $M>0$ such that $|\sigma(y)|\leq M y$, for all $y\in \mathbb{R}$ and
\begin{align*}
|f(x)|&\leq \int_{\mathbb{R}\times \Omega} |a| |\sigma(v^T x + b)| d\pi(a,v,b)\\
&\leq M \int_{\mathbb{R}\times \Omega} |a| |(v^T x + b)| d\pi(a,v,b)\\
&\leq M \int_{\mathbb{R}\times \Omega} |a|(\|x\|_\infty\|v\|_1+|b|) d\pi(a,v,b) 
\end{align*}
for all $x\in X$ and $\pi\in \mathcal{P}(\Omega)$ such that $f=A\pi$. Taking the infimum we get that $|f(x)|\leq M \|x\|_\infty \|f\|_{\mathcal{B}_\sigma} \infty$ for all $f\in \mathcal{B}_\sigma$ and also immediately that the Barron space is an RKBS.

When $\sigma$ is $1$-homogeneous, like for the popular ReLU activation function, the definition of the Barron norm is slightly different
\begin{equation}
\|f\|_{\mathcal{B}_\sigma} := \inf_{f=A\pi} \int_{\mathbb{R}\times \Omega} |a|(\|v\|_1+|b|) d\pi(a,v,b)
\end{equation}

There is a related definition where the above norm is replaced by an infimum over the $L^p(\pi)$-norm of $|a|(1+\|v\|_1+|b|)$ or $|a|(\|v\|_1+|b|)$ for $1$-homogeneous $\sigma$. In this latter case \citet[Proposition 1]{e_barron_2019} shows that the Barron norms for all $1\leq p\leq \infty$ are equal. 

We can show that the Barron space is an instance of the RKBS $\mathcal{F}(X,\Omega)$. First, we treat the case where $\sigma$ is not $1$-homogeneous. 
\begin{theorem}\label{thm:Barron}
Let $\sigma:\mathbb{R}\mapsto \mathbb{R}$ be continuous and let either $\Omega$ be compact or $\sigma(y)\in O(y)$

Let $\mathcal{F}(X,\Omega)$ be the integral RKBS defined by $\varphi(x,w) =  (1+\|v\|_1+|b|)^{-1}\sigma(v^T x+b)$ where $w=(v,b)$. Then $\mathcal{B}_\sigma \cong \mathcal{F}(X,\Omega)$.  
\end{theorem}
\textit{Proof.} Let $\sigma$ be continuous. It is sufficient to show that for every $\pi \in \mathcal{P}(\mathbb{R}\times\Omega)$ there exists a Radon measure $\mu$ such that $A\pi = A\mu$ and 
\[\|\mu\|_{\mathcal{M}(\Omega)}\leq \int_{\mathbb{R}\times \Omega} |a|(1+\|v\|_1+|b|) d\pi(a,v,b)\]
and conversely that for every Radon measure $\mu$ there exists a $\pi\in \mathcal{P}(\mathbb{R}\times\Omega)$ such that $A\mu = A\pi$ and
\[\int_{\mathbb{R}\times \Omega}  |a|(1+\|v\|_1+|b|) d\pi(a,v,b) \leq \|\mu\|_{\mathcal{M}(\Omega)}\]
Taking the infima then gives the equality of norms. In the rest of this proof, we use the convention that $w=(v,b)$.

Define $\varphi(x,w):=(1+\|v\|_1+|b|)^{-1}\sigma(v^T x+b)$. By the assumptions above $\varphi\in C_0(X\times \Omega)$. 

Let $\pi\in \mathcal{P}(\mathbb{R}\times\Omega)$, then define the measure $\mu\in \mathcal{M}(\Omega)$ such that for every Borel set $B\subseteq \Omega$
\[\mu(B) = \int_{\mathbb{R} \times B} a (1+\|v\|_1+|b|) d\pi(a,v,b)\]
Hence, we get that 
\begin{align*}
(A\mu)(x) &= \int_{\Omega} \varphi(x,w) d\mu(w)\\
&= \int_{\mathbb{R}\times \Omega} a (1+\|v\|_1+|b|) (1+\|v\|_1+|b|)^{-1} \sigma(v^Tx+b)d\pi(a,v,b) = (A\pi)(x)    
\end{align*}
for all $x\in X$. Furthermore,
\[|\mu|(\Omega) \leq \int_{\mathbb{R}\times \Omega} |a|(1+\|v\|_1+|b|)d\pi(a,v,b)\]

Let $\mu\in \mathcal{M}(\Omega)$, with its Jordan decomposition $\mu=\mu^+ - \mu^-$. We define the measure $\pi$ as 
\begin{align*}
d\pi(a,v,b) &= \frac{1}{2(1+\|v\|_1+|b|)|\mu^+|(\Omega)} d\mu^+(w) d\delta_{a=2|\mu^+|(\Omega)}(a)\\
&+ \frac{1}{2(1+\|v\|_1+|b|)|\mu^-|(\Omega)} d\mu^-(w) d\delta_{a=2|\mu^-|(\Omega)}(a)
\end{align*}
Clearly, $\pi\in \mathcal{P}(\mathbb{R}\times\Omega)$. Furthermore,
\begin{align*}
    (A\pi)(x) &= \int_{\mathbb{R}\times \Omega} a \sigma(v^T x + b) d\pi(a,v,b)\\
    &= \frac{1}{2|\mu^+|(\Omega)}\int_{\Omega}\int_{\mathbb{R}} a \varphi(x,w)d\delta_{a=2|\mu^+|(\Omega)}(a) d\mu^+(w)\\
    &\quad+ \frac{1}{2|\mu^-|(\Omega)} \int_{\Omega}\int_{\mathbb{R}} a \varphi(x,w) d\delta_{a=2|\mu^-|(\Omega)}(a) d\mu^-(w)\\
    &= \int_{\Omega} \varphi(x,w) d\mu^+(w) + \int_{\Omega} \varphi(x,w) d\mu^-(w)\\
    &= \int_{\Omega} \varphi(x,w) d\mu(w) = (A\mu)(x)
\end{align*}
for all $x\in X$ and 
\begin{align*}
\int_{\mathbb{R}\times \Omega}  |a|(1+\|v\|_1+|b|) d\pi(a,v,b) &= \frac{1}{2|\mu^+|(\Omega)}\int_{\Omega}\int_{\mathbb{R}} |a|d\delta_{a=2|\mu^+|(\Omega)}(a) d\mu^+(w)\\ &+ \frac{1}{2|\mu^-|(\Omega)} \int_{\Omega}\int_{\mathbb{R}} |a| d\delta_{a=2|\mu^-|(\Omega)}(a) d\mu^-(w)\\
&= \mu^+(\Omega) + \mu^-(\Omega) = |\mu|(\Omega)
\end{align*}
Which concludes the proof. \qed

Next, we show that the analogous statement holds when $\sigma$ is $1$-homogeneous. The proof is inspired by \citep[Lemma A.5]{e_kolmogorov_2020} and \citep{e_representation_2021}
\begin{theorem}
Let $\sigma$ be continuous and $1$-homogeneous.  and $\varphi(x,w) = \sigma(v^T x + b)$ with $w=(v,b)$. Then $\mathcal{B}_\sigma\cong \mathcal{F}(X,\Omega)$ 
\end{theorem}
We employ the same technique as the proof of the Theorem~\ref{thm:Barron}. Let $\sigma$ be continuous and $1$-homogeneous, which means that for any $\lambda>0$ and $y\in \mathbb{R}$ it holds that $\sigma(\lambda y) = \lambda \sigma(y)$. The key idea of the proof is reformulating the Barron space over $\Omega'$ where $\|v\|+|b|=1$. 
\[\Omega':= \left.\left\{\frac{w}{\|w\|_1}\right| w\in \Omega, w\neq 0  \right\}\]
Then, as $\Omega'$ is compact, the proof follows from Theorem~\ref{thm:Barron}.

Let $\pi\in \mathcal{P}(\mathbb{R}\times\Omega)$. For any Borel set $B'\subseteq \mathbb{R}\times \Omega'$ we define the Borel set $B \subseteq \mathbb{R}\times \Omega$  
\[B:=\left\{(a,w)\in \mathbb{R}\times \Omega\left| \frac{w}{\|w\|_1}\in B'\right\}\right. \]
We define $\pi'\in \mathcal{P}(\mathbb{R}\times\Omega')$ such that 
\[\pi(B') = \int_{B} (\|v\|_1+|b|) d\pi(a,v,b) \]
Hence, we get by the change of variables formula that for every function $g:\mathbb{R}\times \Omega\mapsto \mathbb{R}$
\[\int_{B} g\left(a,\tfrac{w}{\|w\|_1}\right) (\|v\|_1+|b|) d\pi(a,v,b) = \int_{B'} g'(a,w) d\pi'(a,v,b)\]
where $g'(a,w)= g\left(a,\tfrac{w}{\|w\|_1}\right)$ for all $(a,w)\in \mathbb{R}\times \Omega$. Then we get that
\begin{align*}
(A\pi)(x) &= \int_{\mathbb{R}\times \Omega}a \sigma(v^T x + b) d\pi(a,v,b) \\
&= \int_{\mathbb{R}\times \Omega} a (\|v\|_1+|b|) \sigma\left(\frac{v^Tx+b}{\|v\|_1+|b|}\right)d\pi(a,v,b) \\
&= \int_{\mathbb{R}\times \Omega'} a \sigma(v^Tx+b)d\pi'(a,v,b) \\
&= (A\pi')(x)
\end{align*}
due to the fact that $\sigma$ is $1$-homogeneous. Furthermore,
\begin{align*}
\int_{\mathbb{R}\times \Omega} |a|(\|v\|_1+|b|)d\pi(a,v,b) = \int_{\mathbb{R}\times \Omega'} |a|d\pi'(a,v,b) = \int_{\mathbb{R}\times \Omega'} |a|(\|v\|_1+|b|)d\pi'(a,v,b)    
\end{align*}
 
For the converse, unfortunately we can not say that any $\pi'\in \mathcal{P}(\mathbb{R}\times\Omega')$ is also a probability measure on $\mathbb{R}\times \Omega$, as $\Omega'$ is not necessarily a subset of $\Omega$. Hence, we need to define a function $g:\Omega'\mapsto \mathbb{R}$ such that $g(w)w\in \Omega$, i.e. a scale factor that scales $\Omega'$ into $\Omega$. Then, given a $\pi' \in \mathcal{P}(\mathbb{R}\times\Omega')$, it is easy to see that $\pi$ pushed forward by $g$ is the corresponding measure on $\mathbb{R}\times \Omega$. \qed

\section{Integral RKBS as a Union of p-Norm RKBS}
\label{sec:union}
In this section, we extend Proposition 3 of \citet{e_barron_2021} that the Barron space can be written as a union of RKHS. We show that the integral RKBS $\mathcal{F}(X,\Omega)$ can be written as a union of certain RKBS which are quotient spaces of the $L^p$-spaces. Remarkably, this is independent of the choice $p$. The proof reveals an interesting fact about these RKBS: The curvature of the functions is not relevant for the norm, only the average absolute value of the function.

\subsection{Definitions of p-Norm RKBS}
We first need to properly define the $p$-norm RKBS and the concept of a union of Banach spaces.

Given a $\pi\in \mathcal{P}(\Omega)$, a $1\leq p < \infty$ and a map $\varphi:X\times \Omega \mapsto \mathbb{R}$ such that $\varphi_x\in L^q(\pi)$ for all $x\in X$ and where $q$ is such that $\tfrac{1}{p}+\tfrac{1}{q}=1$. 

We define the $p$-norm RKBS $\mathcal{L}^p(\pi)$ using again the feature map of Definition~\ref{def:rkbs3}. We choose the feature space to be $\Psi=L^p(\pi)$ with its dual space $\Psi^* = L^q(\pi)$, where $\tfrac{1}{p}+\tfrac{1}{q}=1$ and the map $\psi:X\mapsto L^q(\pi)$ as $\psi(x)=\varphi_x$. Effectively, this means that for every $f\in \mathcal{L}^p(\pi)$ there exists an $h\in L^p(\pi)$ such that
\begin{equation}
f(x)=(Ah)(x)=\langle \varphi_x,h\rangle = \int_\Omega \varphi(x,w) h(w) d\pi
\end{equation}
for all $x\in X$. The norm of $f$ is given by the infimum of the $p$-norm of all $h$ which satisfy the above relation. 
\begin{equation}
\|f\|_{\mathcal{L}^p(\pi)} := \inf_{f=Ah} \|h\|_{L^p(\pi)}
\end{equation}

Given an infinite family of normed vector spaces $\mathcal{Z}_I$, indexed by some set $I$, we define the union of these spaces as
\begin{equation}
\mathcal{Z}_{\cup} := \bigcup_{i\in I} \mathcal{Z}_i := \{z\,|\, \exists i \in I \text{ such that }z\in \mathcal{Z}_i\}
\end{equation}
$\mathcal{Z}_{\cup}$ is a normed vector space with the norm 
\begin{equation}
    \|z\|_{\mathcal{Z}_{\cup}} := \inf_{z\in \mathcal{Z}_i} \|z\|_{\mathcal{Z}_i}
\end{equation}

\subsection{Unions of p-Norm RKBS}
We want to show that the RKBS $\mathcal{F}(X,\Omega)$ can be written as a union of $\mathcal{L}^p(\pi)$ for all $\pi\in \mathcal{P}(\Omega)$. We prove this in two steps, first we show the statement holds for $p=1$ and then we show how this extends to any $p>1$.

For this section, we use the assumptions of Section~\ref{sec:notation}. The assumption that $\varphi\in C_0(X\times \Omega)$ implies that $\varphi_x \in L^q(\pi)$ for any $x\in X$, $\pi\in \mathcal{P}(\Omega)$ and $1\leq q\leq \infty$. Hence, for any $1\leq p <\infty$, and $\pi\in \mathcal{P}(\Omega)$ the RKBS $\mathcal{L}^p(\pi)$ is well-defined and so is $\mathcal{F}(X,\Omega)$. 

\begin{theorem}
\begin{equation}
    \mathcal{F}(X,\Omega) \cong \bigcup_{\pi \in \mathcal{P}(\Omega)} \mathcal{L}^1(\pi)
\end{equation}
\end{theorem}
\textit{Proof.} We will prove this theorem by showing that for every Radon measure $\mu \in \mathcal{M}(\Omega)$ there exists a $\pi\in \mathcal{P}(\Omega)$ and function $h\in L^1(\pi)$ such that $Ah = A\mu$ and $\|h\|_{L^1(\pi)}
\leq\|\mu\|_{\mathcal{M}(\Omega)}$ and conversely that for every $\pi\in \mathcal{P}(\Omega)$ and function $h\in L^1(\pi)$ there exists a Radon measure $\mu \in \mathcal{M}(\Omega)$ such that $A\mu = Ah$ and $\|\mu\|_{\mathcal{M}(\Omega)} \leq \|h\|_{L^1(\pi)}$. The equality of the norms then follows from taking the corresponding infima. 

Let $\mu \in \mathcal{M}(\Omega)$ and define the $\pi\in \mathcal{P}(\Omega)$ as
\[\pi := \tfrac{1}{|\mu|(\Omega)} |\mu|\]
As $\mu \ll |\mu| \ll \pi$, we define $h$ as the Radon-Nikodym derivative $h:=\tfrac{d\mu}{d\pi}$. Then 
\[Ah= \int_{\Omega} \varphi(\cdot, w) h(w) d\pi = \int_{\Omega} \varphi(\cdot, w) d\mu = A\mu\]
and 
\[\|h\|_{L^1(\pi)}= \int_{\Omega} |h| d\pi = \int_{\Omega} \text{sgn}(h) h d\pi = \int_{\Omega} \text{sgn}(h) d\mu \leq |\mu|(\Omega)\]
where sgn$(h)$ denotes the sign of $h$.

Conversely, let $\pi\in \mathcal{P}$ and $h\in L^1(\pi)$. Define the positive measures $\mu^+,\mu^-$ such that $d\mu^+ := \max\{h,0\} d\pi$ and $d\mu^- := \min\{h,0\} d\pi$ and define $\mu:=\mu^+-\mu^-$. Clearly, $d\mu = h d\pi$ and
\[A\mu = \int_{\Omega} \varphi(\cdot, w) d\mu = \int_{\Omega} \varphi(\cdot, w) h(w) d\pi = Ah\]
as well as 
\[|\mu|(\Omega) = \int_{\Omega} d\mu^+ + \int_{\Omega} d\mu^- = \int_{\Omega} |h|d\pi = \|h\|_{L^1(\pi)} \] \qed

Next, we show that $\mathcal{L}^1(\pi)$ can be identified as the union of the spaces $\mathcal{L}^p(\tau)$ over all $\tau\in \mathcal{P}(\Omega)$ which are absolutely continuous with respect to $\pi\in \mathcal{P}(\Omega)$. Notably, this result is independent of the value of $p$. This is the case, because for every function $h\in \mathcal{L}^1(\pi)$ we can construct a $\tau\in \mathcal{P}(\Omega)$ such that $h$ becomes constant and for constant functions all $p$-norms are equal.

\begin{theorem}
Let $\pi\in \mathcal{P}(\Omega)$ and let $1\leq p<\infty$, then
\begin{equation}
    \mathcal{L}^1(\pi) \cong \bigcup_{\tau \ll \pi} \mathcal{L}^p(\tau)
\end{equation}
where $\tau \ll \pi$ denotes all $\tau\in \mathcal{P}(\Omega)$ which are absolutely continuous with respect to $\pi$.
\end{theorem}
\textit{Proof.} The proof is similarly structured as the proof of the previous theorem. We will show that for every $h \in L^1(\pi)$ there exists a $\tau\in \mathcal{P}(\Omega)$ and function $h_{p} \in L^p(\tau)$ such that $Ah = Ah_p$ and $\|h_p\|_{L^p(\tau)} \leq \|h\|_{L^1(\pi)}$. Conversely, for every $\tau\in \mathcal{P}(\Omega)$ and function $h_{p} \in L^p(\tau)$ there exists a function $h \in L^1(\pi)$ such that $Ah_p = Ah$ and $\|h\|_{L^1(\pi)} \leq \|h_p\|_{L^p(\tau)}$. Again, taking infima gives the equality of the norms.

Let $h \in L^1(\pi)$.  If $h=0$, then choose $h_p=0$. Otherwise define the measure $\tau$ such that
\[d\tau = \frac{|h|}{\|h\|_{L^1(\pi)}} d\pi\]
Clearly, $\tau \ll \pi$ and a probability measure as
\[\tau(\Omega) = \frac{1}{\|h\|_{L^1(\pi)}} \int_{\Omega} |h(w)| d\pi(w)= 1\]
Next, we define the function
\[h_p(w) = \begin{cases}
 \|h\|_{L^1(\pi)}\frac{h(w)}{|h(w)|} & h(w) \neq 0\\
 0 & h(w) =0
\end{cases}\]
$\tau$ and $g$ are constructed such that $h d\pi = h_p d\tau$, hence
\[Ah= \int_{\Omega} \varphi(\cdot,w) h(w) d\pi(w) = \int_{\Omega} \varphi(\cdot,w) h_p(w) d\tau(w) = Ah_p\]
Furthermore as $|h_p|$ is constant where $h\neq 0$
\[\|h_p\|^p_{L^p(\tau)}= \int_{\Omega} |h_p(w)|^p d\tau(w) = \|h\|^{p-1}_{L^1(\pi)} \int_{\Omega} |h(w)|d\pi(w) = \|h\|^{p}_{L^1(\pi)}\]

Conversely, let $\tau\in \mathcal{P}(\Omega)$ such that $\tau \ll \pi$ and $h_p \in L^p(\tau)$. Define the function $h:= h_p \tfrac{d \tau}{d \pi}$. Then
\[Ah_p = \int_{\Omega} \varphi(\cdot,w) h_p(w) d\tau(w) = \int_{\Omega} \varphi(\cdot,w) h(w) d\pi(w) = Ah \]
and as $\tau$ is a finite measure
\[\|h\|_{L^1(\pi)}= \int_{\Omega} |h(w)|d\pi(w) = \int_{\Omega} |h_p(w)| d\tau(w) = \|h_p\|_{L^1(\tau)} \leq \|h_p\|_{L^p(\tau)}\] \qed

\begin{corollary}
\begin{equation}
    \mathcal{F}(X,\Omega) \cong \bigcup_{\pi \in \mathcal{P}(\Omega)} \mathcal{L}^p(\pi)
\end{equation}
\end{corollary}
Note that due to Theorem~\ref{thm:Barron} and the fact that $\mathcal{L}^2(\pi)$ is an RKHS, we recover the result of \citet{e_barron_2021}.

\section{Duality and Adjoints for RKBS}
\label{sec:duality_rkbs}
\begin{figure}
    \centering
\begin{tikzcd}[column sep=huge,row sep =huge]
\mathcal{M}(\Omega) \arrow[d, "A"']                               & C(\Omega) \arrow[l, "{\langle \mu,g\rangle}", "*"'] \arrow[d, "\subseteq"] \\
{\mathcal{F}(X,\Omega)}                                            & {\mathcal{G}(\Omega,X)} \arrow[l, "{\langle f,g\rangle}", "*"']            \\
C(X) \arrow[r, "{\langle\rho,f\rangle}"', "*"] \arrow[u, "\subseteq"]  & \mathcal{M}(X) \arrow[u, "A^*"']                                     
\end{tikzcd}
    \caption{A diagram depicting the different relations between the Banach spaces of continuous functions $C(X),C(\Omega)$, of Radon measures $\mathcal{M}(X),\mathcal{M}(\Omega)$ and the RKBS $\mathcal{F}(X,\Omega),\mathcal{G}(\Omega,X)$ which are defined in \eqref{eq:F-space} and \eqref{eq:G-space} respectively. The star $*$ denotes a duality relation with the corresponding representation of the duality pairing as used in Lemma \ref{lem:pairing}. The $A$ and $A^*$ denote the embeddings \eqref{eq:A} and \eqref{eq:A*} respectively. The $\subseteq$ relation denotes that one space is a subspace of the other, but the norm is not necessarily inherited.}
    \label{fig:duality_diagram}
\end{figure}
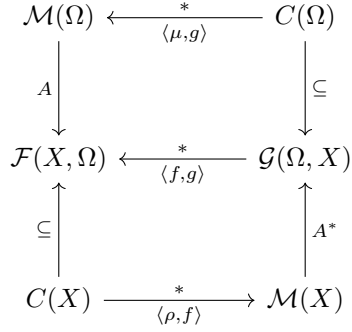

In this section, we characterise the dual space of the class of integral RKBS $\mathcal{F}(X,\Omega)$. Interestingly, the dual space can be interpreted as the continuous functions of the weights which can described by an integral over $\varphi$ with measures over the data instead of the weights. This space is also an RKBS and hence forms an adjoint pair of RKBS together with $\mathcal{F}(X,\Omega)$. Using these adjoints, we formulate a primal-dual optimisation problem for both infinite and finite data. See Figure~\ref{fig:duality_diagram} for an overview of the different spaces and the duality pairings.

\subsection{Duality Pairing}
We define the following Banach space $\mathcal{G}(\Omega,X)$ as a subset of the continuous functions $C_0(\Omega)$ which can be described using Radon measures $\mathcal{M}(X)$ over the data:
\begin{equation}\label{eq:G-space}
\begin{split}
\mathcal{G}(\Omega,X)&:= \{g\in C_0(\Omega) | \exists \rho \in \mathcal{M}(X) \text{ s.t. } A^*\rho = g \}\\
\|g\|_{\mathcal{G}(\Omega,X)} &:= \|g\|_{C_0(\Omega)} 
\end{split}
\end{equation}
Here, the $A^*$ maps data measures $\rho \in \mathcal{M}(X)$ to continuous functions over the weights $C_0(\Omega)$.
\begin{equation}\label{eq:A*}
(A^*\rho)(w) := \int_{X}\varphi(x,w)d\rho(x)\\
\end{equation}
for all $w\in \Omega$ and $\rho \in \mathcal{M}(X)$. The assumption that $\varphi \in C_0(\mathcal{X}\times \Omega)$ in combination with the dominated convergence theorem imply that for all $\rho \in \mathcal{M}(X)$, $A^*\rho \in C_0(\Omega)$ and hence the domain of $A^*$ is the full space of measures $\mathcal{M}(X)$. Also, note that by Definition~\ref{def:rkbs1}, $\mathcal{G}(\Omega,X)$ is an RKBS, as point evaluation is bounded by the supremum norm. 

We will show that $\mathcal{F}(X,\Omega)$ can be identified with the dual space of $\mathcal{G}(\Omega,X)$. The converse is not true, which is a consequence of the fact that the space of continuous functions $C_0(\Omega)$ is not a reflexive. First, we will show that there exists a well-defined pairing between $\mathcal{G}(\Omega,X)$ and $\mathcal{F}(X,\Omega)$.

\begin{lemma}[Data Weight Pairing]\label{lem:pairing}
Consider the bi-linear map defined as
\begin{equation} \label{eq:duality_pairing}
    \langle f,g \rangle := \langle \rho, f \rangle = \langle \mu, g \rangle = \langle \rho\times \mu, \varphi \rangle = \int_{X\times \Omega} \varphi(x,w) d(\rho \times \mu)(x,w) 
\end{equation}
for all $g\in \mathcal{G}(\Omega,X), f\in \mathcal{F}(X,\Omega)$ and where $\rho\in \mathcal{M}(X),\mu\in \mathcal{M}(\Omega)$ such that $g=A^*\rho,f=A\mu$. Here $\langle \rho, f \rangle$ is a pairing between $C_0(X)$ and $\mathcal{M}(X)$, $\langle \mu, g \rangle$ between $C_0(\Omega)$ and $\mathcal{M}(\Omega)$ and $\langle \rho\times \mu, \varphi \rangle$ between $C_0(X\times \Omega)$ and $\mathcal{M}(X\times \Omega)$. (See also Figure~\ref{fig:duality_diagram})

This map is well-defined and bounded
\begin{equation}
|\langle f,g \rangle| \leq \|f\|_{\mathcal{F}(X,\Omega)}\|g\|_{\mathcal{G}(\Omega,X)}
\end{equation}
\end{lemma}

\textit{Proof.} First, we show that \eqref{eq:duality_pairing} is well-defined. Let $g\in \mathcal{G}(\Omega,X), f\in \mathcal{F}(X,\Omega)$ and let $\rho\in \mathcal{M}(X),\mu\in \mathcal{M}(\Omega)$ such that $g=A^*\rho,f=A\mu$. By the assumptions in Section~\ref{sec:notation}, $\varphi \in C_0(X\times \Omega)$ and the product of two Radon measures is again a Radon measure, hence 
\[\int_{X\times \Omega} \varphi(x,w) d(\rho \times \mu)(x,w) = \langle \rho\times \mu, \varphi \rangle\]
is well-defined and finite, where the pairing is between $C_0(X\times \Omega)$ and $\mathcal{M}(X\times \Omega)$. This means we can use Fubini's theorem
\[\int_{X\times \Omega} \varphi(x,w) d(\rho \times \mu)(x,w)  = \int_X \int_\Omega \varphi(x,w) d\mu(w) d\rho(x) = \int_X f(x) d\rho(x) = \langle \rho, f \rangle \]
where this last pairing is between $C_0(X)$ and $\mathcal{M}(X)$, as well as 
\[\int_{X\times \Omega} \varphi(x,w) d(\rho \times \mu)(x,w) = \int_\Omega \int_X \varphi(x,w) d\rho(x) d\mu(w) = \int_\Omega g(w) d\mu(w) = \langle \mu, g \rangle \]
where this last pairing is between $C_0(\Omega)$ and $\mathcal{M}(\Omega)$.

We can show that $\langle f,g \rangle$ is independent on the choice of $\mu$ and $\rho$. If $\rho'\in \mathcal{M}(X),\mu'\in \mathcal{M}(\Omega)$ such that $g=A^*\rho',f=A\mu'$ then 
\begin{align*}
\int_{X\times \Omega} \varphi(x,w) d(\rho' \times \mu')(x,w) &= \langle \rho', f \rangle = \int_{X\times \Omega} \varphi(x,w) d(\rho' \times \mu)(x,w)\\
&= \langle \mu, g \rangle = \int_{X\times \Omega} \varphi(x,w) d(\rho \times \mu)(x,w)    
\end{align*}
Hence, the map is well-defined. Furthermore, the pairing is bounded as
\[|\langle f,g \rangle| = |\langle \mu, g \rangle| \leq \|\mu\|_{\mathcal{M}(\Omega)} \|g\|_{C_0(\Omega)}\]
for all $\mu \in \mathcal{M}(\Omega)$ such that $f=A\mu$. Hence, by taking the corresponding infimum we get the inequality
\[|\langle f,g \rangle| \leq \|f\|_{\mathcal{F}(X,\Omega)}\|g\|_{\mathcal{G}(\Omega,X)}\]
\qed

This lemma also implies that $A^*$ can be understood as an adjoint of $A$ as for all $\rho\in \mathcal{M}(X),\mu\in \mathcal{M}(\Omega)$
\begin{equation}
    \langle \rho, A\mu \rangle = \langle \mu, A^* \rho \rangle = \langle A^*\rho, A\mu \rangle
\end{equation}
Here the first pairing is between $C_0(X)$ and $\mathcal{M}(X)$, the second between $C_0(\Omega)$ and $\mathcal{M}(\Omega)$ and the third between $\mathcal{G}(\Omega,X)$ and $\mathcal{F}(X,\Omega)$.

With this pairing, we can show that $\mathcal{F}(X,\Omega)$ can indeed be identified as the dual of $\mathcal{G}(\Omega, X)$. We prove this by leveraging the fact that the dual of a subspace can be identified as the quotient space with respect to its annihilator. As $\mathcal{F}(X,\Omega)$ is isomorphic to the quotient space $\mathcal{M}(\Omega)/\mathcal{N}(A)$ by definition, the proof consists of showing that $\mathcal{N}(A)$ is indeed the annihilator of $\mathcal{G}(\Omega, X)$.

\begin{theorem}
Let $\mathcal{G}^*(\Omega,X)$ be the dual space of $\mathcal{G}(\Omega,X)$, then $\mathcal{G}^*(\Omega,X) \cong \mathcal{F}(X,\Omega)$
\end{theorem}
\textit{Proof.}
The dual of a normed subspace $\mathcal{Y}$ of a Banach space $\mathcal{Z}$ can be identified as a quotient space of the dual space of the full space $\mathcal{Z}^*$ over the annihilator of the subspace $\mathcal{Y}^\perp$ \citep[Theorem 4.8]{rudin_functional_1991}. Here, the annihilator $\mathcal{Y}^\perp$ is defined as
\begin{equation}
    \mathcal{Y}^\perp := \{z^* \in \mathcal{Z}^* \,|\, \forall y\in \mathcal{Y}, \langle z^*,y\rangle =0 \}
\end{equation}
where the pairing is between the spaces $\mathcal{Z}$ and $\mathcal{Z}^*$. Remember, that $\mathcal{F}(X,\Omega)$ is defined as the quotient space $\mathcal{M}(\Omega)/\mathcal{N}(A)$. So it only remains to prove that $\mathcal{G}(\Omega,X)^\perp \cong \mathcal{N}(A)$, where
\begin{equation}
    \mathcal{G}(\Omega,X)^\perp := \{\mu \in \mathcal{M}(\Omega) \,|\, \forall g \in \mathcal{G}(\Omega, X), \langle \mu, g \rangle =0 \}
\end{equation}
The isometry is trivial as both spaces use the total variation norm of $\mathcal{M}(\Omega)$. We prove that these spaces are isomorphic in two steps.  

Let $\mu \in \mathcal{G}(\Omega,X)^\perp$, $\rho \in \mathcal{M}(X)$ and define the function $g:=A^* \rho \in \mathcal{G}(\Omega,X)$. Using Lemma~\ref{lem:pairing}
\[0=\langle \mu, g\rangle = \langle A\mu, \rho \rangle \]
As this holds for all $\rho\in \mathcal{M}(X)$ we get that $A\mu\equiv 0$ and hence $\mu \in \mathcal{N}(A)$.

Conversely, let $\mu \in \mathcal{N}(A)$ and $g \in \mathcal{G}(\Omega, X)$. Furthermore, let $\rho\in \mathcal{M}(X)$ such that $A^*\rho = g$. Again by Lemma~\ref{lem:pairing}
\[0=\langle A\mu, \rho \rangle = \langle \mu, g \rangle \]
Hence, $\mu \in \mathcal{G}(\Omega,X)^\perp$ \qed

As both $\mathcal{F}(X,\Omega)$ and $\mathcal{G}(\Omega, X)$ are RKBSs, we can show that they form an adjoint pair of RKBSs and $\varphi$ is the reproducing kernel in the sense of Definition~\ref{def:rkbs2}, as $\varphi(x,w) = \langle A\delta_w, A^*\delta_x \rangle$.

\begin{corollary}
$\mathcal{F}(X,\Omega)$ and $\mathcal{G}(\Omega, X)$ form an adjoint pair of RKBS with the reproducing kernel $K=\varphi$
\end{corollary}
\textit{Proof.}
Define $K:=\varphi$. Then 
\[K(x,\cdot) = \varphi_x = \int_X \varphi(x',\cdot) d\delta_x(x') = A^*\delta_x \in \mathcal{G}(\Omega, X)\]
as $\delta_x\in \mathcal{M}(X)$ for all $x\in X$, i.e. all point measures are Radon measures. Furthermore, let $f\in \mathcal{F}(X,\Omega)$, then by Lemma~\ref{lem:pairing}
\[\langle f, K(x,\cdot) \rangle = \langle f, A^*\delta_x \rangle = \langle \delta_x, f \rangle = f(x) \]
for all $x\in X$.

For the adjoint we find that
\[K(\cdot, w) = \varphi_w = \int_\Omega \varphi(\cdot,w') d\delta_w(w') = A\delta_w \in \mathcal{F}(X,\Omega)\]
as $\delta_w \in \mathcal{M}(\Omega)$ for all $w\in \Omega$. Furthermore, let $g\in \mathcal{G}(\Omega,X)$ then again by Lemma~\ref{lem:pairing}
\[\langle K(\cdot,w),g \rangle = \langle A\delta_w , g \rangle = \langle \delta_w,g\rangle = g(w)\]
for all $w\in \Omega$. \qed

\subsection{Dual Formulation of ERM with Infinite Data}
The description of the dual space allows a dual formulation of the corresponding empirical risk optimisation (ERM) problem for integral RKBS. This opens the door to many primal-dual type algorithms to be used for these problems. First, we will consider the infinite data case, before sampling to finite data. 

Given data pairs $(x,y(x))$ for all $x\in X$ and a probability measure $\nu \in \mathcal{M}(X)$ where it is assumed that $y\in L^2(\nu)$, we define the ERM as 
\begin{equation}\label{eq:ERM}
    \inf_{\mu \in \mathcal{M}(\Omega)} J(A\mu) + R(\mu)
\end{equation}
where the data fidelity term $J:\mathcal{F}(X,\Omega)\mapsto \mathbb{R}$ and the regularisation term $R:\mathcal{M}(\Omega)$ are defined by
\begin{equation}
\begin{split}
J(A\mu) &:= \tfrac{1}{2}\|A\mu - y\|_{L^2(\nu)}^2 = \frac{1}{2} \int_{X} ((A\mu)(x)-y(x))^2 d\nu(x) \\
R(\mu) &:= \|\mu\|_{\mathcal{M}(\Omega)}= |\mu|(\Omega)
\end{split}
\end{equation}

Note that for the infinite setting, the optimisation problem is convex. The dual problem can be constructed in terms of the convex conjugates of $J$ and $R$, also referred to as the Fenchel–Legendre transform. However, as $\mathcal{M}(\Omega)$ and $\mathcal{F}(X,\Omega)$ are already dual spaces of a non-reflexive Banach space, we define the convex conjugates $J^*:\mathcal{M}(X)\mapsto \mathbb{R}$ and $R^*:\mathcal{G}(\Omega,X)\mapsto \mathbb{R}$ in terms of the pre-dual
\begin{equation}
\begin{split}
J^*(\rho) &:= \sup_{f \in \mathcal{F}(X,\Omega)}(\langle \rho, f \rangle - J(f))\\
R^*(g) &:= \sup_{\mu \in \mathcal{M}(\Omega)}(\langle \mu, g \rangle - R(\mu))
\end{split}
\end{equation}
where $\rho \in \mathcal{M}(X)$ and $g\in C_0(\Omega)$. We can write $J^*$ as a function of $\rho\in \mathcal{M}(X)$ instead of $g\in \mathcal{G}(\Omega,X)$ as $\langle f, A^*\rho \rangle=\langle \rho, f \rangle$ by Lemma~\ref{lem:pairing}.

\begin{lemma}\label{lem:dual_problem}
If $\mathcal{F}(X,\Omega)$ is dense in $L^2(\nu)$, then the convex conjugates of $J$ and $R$ are given by
\begin{equation}
\begin{split}
J^*(\rho) &= \begin{cases} \int_X \frac{1}{2}\frac{d\rho}{d\nu}(x) + y(x) d\rho(x) & \rho \ll \nu \\ \infty & \text{otherwise} \end{cases}\\
R^*(g) &= \begin{cases} 0 & \|g\|_{C_0(\Omega)}\leq 1\\ \infty & \text{otherwise} \end{cases}
\end{split}
\end{equation}
for all $\rho\in \mathcal{M}(X)$ and $g\in C_0(X)$
\end{lemma}
\textit{Proof.}
First, we prove the identity for $J^*$. Let $\rho\in \mathcal{M}(X)$. We consider $\rho \ll \nu$. Then for some $f\in\mathcal{F}(X,\Omega)$
\begin{align*}
\langle \rho, f \rangle - J(f) &= \int_X f(x) d\rho(x) -\frac{1}{2} \int_{X} (f(x)-y(x))^2 d\nu(x)\\
&=\int_X f(x) \frac{d\rho}{d\nu}(x) -\frac{1}{2} (f(x)-y(x))^2 d\nu(x)\\
&= \int_X -\frac{1}{2} \left(f(x)-y(x)-\frac{d\rho}{d\nu}(x)\right)^2 + \frac{d\rho}{d\nu}(x)\left(\frac{1}{2}\frac{d\rho}{d\nu}(x) + y(x)\right)  d\nu(x)\\
\end{align*}
where we have completed the square. By the assumption that $\mathcal{F}(X,\Omega)$ is dense $L^2(\nu)$, the square vanishes in the supremum, and we get
\[\int_X \frac{d\rho}{d\nu}(x)\left(\frac{1}{2}\frac{d\rho}{d\nu}(x) + y(x)\right) d\nu(x) = \int_X \frac{1}{2}\frac{d\rho}{d\nu}(x) + y(x) d\rho(x)\]
Next, when $\rho$ is not absolutely continuous with respect to $\nu$, then there exists a Borel set $B$ of $X$ such that $\nu(B)=0$ but $\rho(B)\neq 0$. Then for a sequence of $f_{n}\in\mathcal{F}(X,\Omega)$, $n\in \mathbb{N}$, converging to the indicator function of $B$
\[\lim_{c\mapsto \infty} \lim_{n\mapsto \infty} (\langle \rho, c f_n \rangle - J(cf_n)) = \lim_{c\mapsto \infty} c\rho(B) = \infty\]

Finally, we prove  the identity for$R^*$, which is a standard argument. Let $g\in C_0(\Omega)$. First, we consider the case when $g$ is in the unit ball, i.e. $\|g\|_{C_0(\Omega)}\leq 1$. Then for any $\mu \in \mathcal{M}(\Omega)$ we find that
\[\langle \mu, g \rangle - R(\mu) \leq (\|g\|_{C_0(\Omega)}-1)\|\mu\|_{\mathcal{M}(\Omega)}\leq 0\] 
Hence, the supremum is achieved when $\mu$ is the zero measure. Next, when $\|g\|_{C_0(\Omega)}> 1$, there exists a $w\in \Omega$ such that $|g(w)|>1$. Taking the point measures $c \tfrac{|g(w)|}{g(w)} \delta_w$ then 
\[\lim_{c\mapsto \infty}\left( \left\langle c \tfrac{|g(w)|}{g(w)} \delta_w , g \right\rangle - R\left(c \tfrac{|g(w)|}{g(w)} \delta_w\right) \right) = \lim_{c\mapsto \infty} c(|g(w)|-1|) = \infty \] \qed

By Fenchel's duality theorem, we get the strong duality equality for \eqref{eq:ERM}. 

\begin{theorem}\label{thm:strong_duality}
The strong duality equality holds
\begin{equation}
\inf_{\mu \in \mathcal{M}(\Omega)} J(A\mu) + R(\mu) = \sup_{\rho \in \mathcal{M}(X)} - J^*(-\rho) - R^*(A^*\rho)    
\end{equation}
The perturbation function for this  problem is given by 
\begin{equation}
F(\mu,f) = J(A\mu - f) + R(\mu)
\end{equation}
where $\mu\in \mathcal{M}(\Omega)$ and $f\in \mathcal{F}(X,\Omega)$.
\end{theorem}
\textit{Proof.} By Fenchel's duality theorem, strong duality holds if there exists a $\mu \in \mathcal{M}(\Omega)$ such that $R(\mu)$ is finite and $J$ is finite and continuous at $A\mu$. Clearly, $R(\mu)$ is finite for any $\mu\in \mathcal{M}(\Omega)$ and $J$ is continuous by dominated convergence as $A\mu \in C_0(X)$. \qed

\subsection{Dual Formulation of ERM with Finite Data}
If we have finite data, i.e. $X=\{x_1,\cdots,x_n\}$, then we can rewrite the dual problem in simpler terms. When $X$ is finite, then the dimension of $\mathcal{M}(X)$ is also finite and consequently so are $\mathcal{G}(\Omega,X)$ and $\mathcal{F}(X,\Omega)$. Due to the fact that the point measures are the extreme points of the unit ball of $\mathcal{M}(\Omega)$, the infimum in the definition of the norm of $\mathcal{F}(X,\Omega)$ \eqref{eq:norm_f} has a minimiser which is a linear combination of point measures $\delta_{w_i}$, for certain $w_i\in \Omega$ and $i\in \{1,\cdots, n\}$. This is the basis for the representer theorems of \citet{bartolucci_understanding_2023,bredies_sparsity_2019}
\begin{theorem}\citep[Theorem 3.9]{bartolucci_understanding_2023}\label{thm:representer}
Let $X=\{x_1,\cdots,x_n\}$, for some $n\in \mathbb{N}$. Then the ERM \eqref{eq:ERM} admits a minimiser $\mu^*$ of the form
\begin{equation}
\mu^* = \sum_{i=1}^n \alpha_i \delta_{w_i}
\end{equation}
for some $\alpha_i\in \mathbb{R}, w_i\in \Omega$ where $i\in \{1,\cdots, n\}$ 
\end{theorem}

Note that $f^*:=A\mu^*$ can be written in terms of the reproducing kernel $K=\varphi$
\[f^*(x) := (A\mu^*)(x) = \sum_{i=1}^n \alpha_i \varphi(x,w_i) \]
The trade-off for the representer theorem is that finding the $w_i$ is a non-convex problem. So we have traded the convexity of the ERM for finite dimensionality, i.e. we only need to find $2n$ parameters. This is specific to RKBS, in contrast to RKHS, where the finite dimensional problem remains convex. 

As $X$ is finite, we have finite data pairs $(x_i,y_i)$ and the probability measure $\nu := \sum_{i=1}^n \gamma_i \delta_{x_i}$, where $\gamma_i>0$ for $i\in \{1,\cdots, n\}$ and $\sum_{i=1}^n \gamma_i=1$. Similarly, we can write any measure $\rho \in \mathcal{M}(X)$ as $\rho := \sum_{i=1}^n \lambda_i \delta_{x_i}$, where $\lambda_i\in \mathbb{R}$ for $i\in \{1,\cdots, n\}$. The dual ERM problem then becomes
\begin{theorem}
Let $X$ have cardinality $n<\infty$, then we can rewrite the primal ERM \eqref{eq:ERM} using
\begin{equation}
\begin{split}
J(\alpha,\mathbf{w}) &= \sum_{i,j=1}^n \gamma_i (\alpha_j \varphi(x_i, w_j) - y_i)^2\\
R(\alpha) &= \|\alpha\|_1    
\end{split}
\end{equation}
for $\alpha=(\alpha_1,\cdots, \alpha_n)\in \mathbb{R}^n, \mathbf{w}= (w_1,\cdots, w_n)\in \Omega^n$.

The dual problem can be written using
\begin{equation}
\begin{split}
J^*(\lambda) &= \sum_{i=1}^n \lambda_i\left(\frac{\lambda_i}{2 \gamma_i} + y_i\right) \\
R^*(\lambda) &= \begin{cases} 0 & -1 \leq \sum_{i=1}^n \lambda_i \varphi(x_i,w) \leq 1 \quad \forall w\in \Omega \\ \infty & \text{otherwise} \end{cases}
\end{split}
\end{equation}
for $\lambda=(\alpha_1,\cdots, \alpha_n)\in \mathbb{R}^n$.

The strong duality equality holds
\[\inf_{\alpha\in \mathbb{R}^n, \mathbf{w}\in \Omega^n} J(\alpha, \mathbf{w}) + R(\alpha) = \sup_{\lambda \in \mathbb{R}^n} - J^*(-\lambda) - R^*(\lambda) \]
\end{theorem}
\textit{Proof.} The statement follows directly from Lemma~\ref{lem:dual_problem} and Theorems~\ref{thm:strong_duality} and~\ref{thm:representer}. Note that by the assumption that $\gamma_i>0$ for $i\in \{1,\cdots,n\}$, all measures $\rho \in \mathcal{M}(X)$ are absolutely continuous with respect to $\nu$ and
\[\frac{d\rho}{d\nu}(x_i) = \frac{\lambda_i}{\gamma_i}\]
for every $i\in \{1,\cdots,n\}$. \qed

\section{Conclusion and Discussion}
In this paper, we have developed duality theory for neural networks, including the Barron spaces of \citet{e_priori_2019} and the integral RKBS of \citet{bartolucci_understanding_2023}. A key conceptual result is that the weights and the data are dual concepts. On one side, we have neural networks which represent functions of the data parameterised by some weights, where on the other side we have the dual networks which are functions of the weights parameterised by some data points. These dual networks can be thought of as representing the range of functions which can be represented given some sample of the data, which is closely related to complexity measures, such as the Rademacher complexity. 

We have considered RKBS instead of more commonly used RKHS. We have shown that spaces for neural networks such as the Barron spaces can be properly described by RKBS. However, a common problem which involves RKBS is that, when considering finite data samples, the ERM optimisation problem becomes non-convex. One of the main advantages of RKHS-based methods, such as support vector machines, is that the ERM is convex and thus relatively easy to solve. This is in contrast to neural networks, where gradient descent based methods are used. The class of integral RKBS we have investigated, however, is a much more flexible function space, as it is a union of a large class of different RKHS.

We have shown that the Barron space $\mathcal{B}_\sigma$ can be identified with the integral RKBS $\mathcal{F}(X,\Omega)$. This theorem holds for the most commonly used activation functions $\sigma$, like the ReLU, sigmoid, inverse tangent and sine activation functions are valid. For activation functions which are higher powers of the ReLU function, special considerations are necessary. However, with the proper definition of the Barron norm for these activation functions, we conjecture that a similar statement holds. 

The dual viewpoint has multiple use cases. It can be used to bring a foundation to architecture search methods for neural networks \citep{bungert_neural_2021}. These methods are based on Bregman iteration, which are closely related to augmented Lagrangian techniques \citep{brune_primal_2011}. Another use case is related to sampling. The question of which data points to sample is often referred to as active learning or sampling in experimental design \citep{settles_active_2009}. The functions in the dual space enable a way to address this problem for neural networks. On the other hand, adversarial attacks try to find input data which gives an undesired result \citep{madry_towards_2019}. In the language of this paper, this means, given slightly different data measures, leads to significant differences in the dual space. Therefore, this work can give a new viewpoint on quantifying adversarial robustness.   

One of the key advances in the field of neural networks was the change from shallow to deep neural networks. Most literature dealing with the theoretical study of neural networks has so far been restricted to shallow neural networks \citep{e_priori_2019,parhi_banach_2021}. Recently, \citet{unser_representer_2019} proved a representer theorem for certain deep networks with the ReLU activation function and \citet{e_banach_2020} developed generalised Barron spaces which deal with deep neural networks. In principle, there are no barriers to considering deep neural networks in our integral RKBS framework \eqref{eq:F-space}, by taking the kernel $\varphi$ to be a concatenation of activation functions $\sigma$ and affine functions based on parameters $w^j$ for layers $j\in\{1,\cdots,L\}$. Then we can still take measures $\mu \in \mathcal{M}(\Omega)$ where the parameters are from $\Omega=\Omega_1 \times \cdots \times \Omega_L$. We are not sure if the resulting norm of $\mathcal{F}(X,\Omega)$ is appropriate for deep neural networks, as it is then mostly dependent on the size of weights in the last (linear) layer. However, we conjecture that the ideas of the proof~\ref{thm:Barron} can be extended to prove a similar result for generalised Barron spaces. 

\section*{Acknowledgements}
C.B. acknowledges support by the European Union's Horizon 2020 research and innovation programme under the Marie Skłodowska-Curie grant agreement No 777826 (NoMADS). L.S. thanks Stephan van Gils for his support and encouragement for taking the research in this direction.

\bibliographystyle{plainnat}
\bibliography{NeuralNetwork} 

\end{document}